\documentclass[12pt,english]{article}
\usepackage{theorem}
\usepackage{babel}
\usepackage{a4wide}
\usepackage{graphics}
\usepackage{epsfig}
\usepackage{amssymb}
\usepackage{latexsym}
\usepackage{amsmath}

\newtheorem{thh}{Theorem}[section]
\newtheorem{df}[thh]{Definition}

\newtheorem{lem}[thh]{Lemma}
\newtheorem{cor}[thh]{Corollary}
\newtheorem{prop}[thh]{Proposition}

\title{Relative exactness modulo a polynomial map and algebraic
$\CP$-actions}
\author{Philippe Bonnet}

\newcommand{\dem}{{\em Proof: }}
\newcommand{\qed}{\begin{flushright} $\blacksquare$\end{flushright}}
\newcommand{\der}{ \partial}
\newcommand{\D}{{\cal{D}}}
\newcommand{\CC}{\mathbb C}
\newcommand{\TF}{{\cal{T}}^1(F)}

\newcommand{\OF}{\omega_F}
\newcommand{\CX}{ \mathbb C [x_1,...,x_n]}

\newcommand{\CF}{ \mathbb C [F]}
\newcommand{\CT}{ \mathbb C [t_1,...,t_q]}
\newcommand{\CN}{\mathbb C ^n}
\newcommand{\CQ}{\mathbb C ^q}
\newcommand{\CP}{(\mathbb C ^p,+)}

\begin{document}
\maketitle

\begin{center}
Section de Math\'ematiques, Universit\'e de Gen\`eve, \\
2-4, rue du li\`evre,
1211 Gen\`eve 24, Switzerland \\ e-mail: Philippe.bonnet@math.unige.ch

\end{center}

\begin{abstract}
Let $F=(f_1,..,f_q)$ be a polynomial dominating map from $\CN$ to
$\CQ$. In this paper we study the quotient $\TF$ of polynomial
1-forms that are exact along the generic fibres of $F$, by 1-forms
of type $dR + \sum a_i df_i$, where $R,a_1,..,a_q$ are
polynomials. We prove that $\TF$ is always a torsion $\CT$-module.
Then we determine under which conditions on $F$ we have $\TF=0$.
As an application, we study the behaviour of a class of algebraic
$\CP$-actions on $\CN$, and determine in particular when these
actions are trivial.
\end{abstract}

\section{Introduction}

Let $F=(f_1,..,f_q)$ be a dominating polynomial map from $\CN$ to
$\CQ$ with $n>q$. Let $\Omega^k(\CN)$ be the space of polynomial
differential $k$-forms on $\CN$. For simplicity, we denote by
$\CC[F]$ the algebra generated by $f_1,..,f_q$, and by $\CC(F)$
its fraction field. Our purpose in this paper is to compare two
notions of relative exactness modulo $F$ for polynomial 1-forms,
and to deduce some consequences on some algebraic groups actions.

The first notion is the {\em topological relative exactness}. A
polynomial 1-form $\omega$ is topologically relatively exact (in
short: TR-exact) if $\omega$ is exact along the generic fibres of
$F$. More precisely this means there exists a Zariski open set $U$
in $\CQ$ such that, for any $y$ in $U$, the fibre $F^{-1}(y)$ is
non-critical and non-empty, and $\omega$ has null integral along
any loop $\gamma$ contained in $F^{-1}(y)$.

The second notion is the {\em algebraic relative exactness}. A
polynomial 1-form is algebraically relatively exact (in short:
AR-exact) if it is a coboundary of the De Rham relative complex of
$F$ (\cite{Ma2}). Recall this complex is given by the spaces of
relative forms: $$ \Omega ^k _F= \Omega ^k(\CN)/\sum df_i \wedge
\Omega^{k-1}(\CN) $$ and the morphisms $d_F:\Omega ^k
_F\longrightarrow \Omega ^{k+1} _F$ induced by the exterior
derivative.

\begin{df}
The module of relative exactness of $F$ is the quotient $\TF$ of
TR-exact 1-forms by AR-exact 1-forms. This is a $\CC[F]$-module
under the multiplication rule $(P(F),\omega)\mapsto P(F)\omega$.
\end{df}
 For holomorphic germs, Malgrange implicitly
compared these notions of relative exactness in \cite{Ma2}. He
proved that the first relative cohomology group of the germ $F$ is
zero if the singular set of $F$ has codimension $\geq 3$; in this
case, ${\cal{T}}^1(F)$ is reduced to zero. In \cite{B-C}, Berthier
and Cerveau studied the relative exactness of holomorphic
foliations, and introduced a similar quotient. For polynomials in
two variables, Gavrilov proved that ${\cal{T}}^1(f)=0$ if every
fibre of $f$ is connected and reduced (\cite{Ga}). Concerning
polynomial maps, we first prove the following result.

\begin{prop} \label{torsion}
If $F$ is a dominating map, then $\TF$ is a torsion $\CF$-module.
\end{prop}
In other words, every TR-exact 1-form $\omega$ can be written as:
$$P(F)\omega = dR +a_1df_1 + ..+ a_qdf_q$$ where $R,a_1,..,a_q$
are all polynomials. In \cite{B-D}, the author in collaboration
with Alexandru Dimca studied in a comprehensive way the torsion of
this module for any polynomial function $f:\CC^2 \rightarrow \CC$.
We are going to extend these results in any dimension and
determine when $\TF$ is zero.

Let $F:X\rightarrow Y$ be a morphism of algebraic varieties, where
$Y$ is equidimensionnal and $X$ may be reducible. A property
${\cal{P}}$ on the fibres of $F$ is {\em $k$-generic} if the set
of points $y$ in $Y$ whose fibre $F^{-1}(y)$ does not satisfy
${\cal{P}}$ has codimension $>k$ in $Y$. A {\em blowing-down} is
an irreducible hypersurface $V$ in $\CN$ such that $F(V)$ has
codimension $\geq 2$ in $\CQ$. If no such hypersurface exists, we
say that $F$ has no blowing-downs. Finally $F$ is non-singular in
codimension 1 if its singular set has codimension $\geq 2$.  It is
easy to prove that a non-singular map in codimension 1 has no
blowing-downs.

\begin{df}
The map $F$ is primitive if its fibres are 0-generically connected
and 1-generically non-empty.
\end{df}
Then we show that a polynomial map $F$ is primitive if and only if
every polynomial $R$ locally constant along the generic fibres of
$F$ can be written as $R=S(F)$, where $S$ is a polynomial. So this
definition extends the notion of primitive polynomial
(\cite{D-P}).

\begin{df} \label{quasi}
The map $F$ is quasi-fibered if $F$ is non-singular in codimension 1,
its fibres are 1-generically connected and 2-generically non-empty.
The map $F$ is weakly quasi-fibered if $F$ has no blowing-downs, its
fibres are 1-generically connected and 2-generically non-empty.
\end{df}

\begin{thh} \label{quasi3}
Let $F$ be a primitive mapping. If $F$ is a quasi-fibered mapping,
then $\TF=0$. If $F$ is weakly quasi-fibered, then every TR-exact
1-form $\omega$ splits as $\omega=dR + \omega_0$, where $R$ is a
polynomial and $\omega_0 \wedge df_1 \wedge ..\wedge df_q=0$.
\end{thh}
We apply these results to the study of algebraic $\CP$-actions on
$\CN$. Such an action is a regular map $\varphi: \CC^p \times \CN
\rightarrow \CN$ such that
$\varphi(u,\varphi(v,x))=\varphi(u+v,x)$ for all $u,v,x$.
Geometrically speaking, $\varphi$ is obtained by integrating a
system $\D=\{\der_1,..,\der_p\}$ of derivations on $\CX$ that are
pairwise commuting and locally nilpotent (\cite{Kr}), that is :
$$\forall f \in \CX, \quad \exists k \in \mathbb{N}, \quad  \der_i
^k (f)=0$$ The ring of invariants $\CX ^{\varphi}$ is the set of
polynomials $P$ such that $P\circ \varphi=P$. Finally $\varphi$ is
{\em free at the point $x$} if the orbit of $x$ has dimension $p$,
and {\em free} if it is free at any point of $\CN$. The set of
points where $\varphi$ is not free is an algebraic set denoted
${\cal{NL}}(\varphi)$.

\begin{df}
An algebraic $\CP$-action on $\CN$ satisfies condition $(H)$ if
its ring of invariants is isomorphic to a polynomial ring in
$(n-p)$ variables.
\end{df}
Under this condition, $\varphi$ is provided with a {\em quotient
map} $F$ (\cite{Mu2}) defined as follows: If $f_1,.. ,f_{n-p}$
denote a set of generators of $\CX ^{\varphi}$, then: $$ F: \CN
\longrightarrow \CC^{n-p},\; x\longmapsto (f_1(x),..,f_{n-p}(x))
$$ The generic fibres of $F$ are orbits of the action, but this
map need not define a topological quotient: For instance, it does
not separate all the orbits. The action $\varphi$ is {\em trivial}
if it is conjugate by a polynomial automorphism of $\CN$ to the
action: $$ \varphi_0(t_1,..,t_p;x_1,..,x_n)= (x_1 + t_1,..,x_p
+t_p,x_{p+1},..,x_n) $$ We are going to search under which
conditions the actions satisfying $(H)$ are trivial. According to
a result of Rentschler (\cite{Re}), every fix-point free algebraic
$(\CC,+)$-action on $\CC^2$ is trivial. We know that $(H)$ is
always satisfied for $(\CC,+)$-actions on $\CC^3$ (\cite{Miy}),
but we still do not know if fixed-point free $(\CC,+)$-actions on
$\CC^3$ are trivial (\cite{Kr}). In dimension $\geq 4$, the works
of Nagata and Winkelmann (\cite{Kr},\cite{Wi}) prove that $(H)$
need not be satisfied. For $(\CC,+)$-actions satisfying this
condition, Deveney and Finston proved that $\varphi$ is trivial if
its quotient map defines a locally trivial $(\CC,+)$-fibre bundle
on its image ([D-F]).

We are going to see how this last result extends via relative
exactness. Let $\varphi$ be a $\CP$-action on $\CN$ satisfying
$(H)$, and consider the following operators: $$
\begin{array}{rcl}
[\D]: & (R_1,..,R_p) \longmapsto & \det((\partial_i(R_j))) \\ \\

J: & (R_1,..,R_p)\longmapsto & \det(dR_1,..,dR_p,df_1,..,df_{n-p})
\end{array}
$$ We say that $[\D]$ (resp. $J$) vanishes at the point $x$ if, for any
polynomials $R_1,..,R_p$, we have $[\D](R_1,..,R_p)(x)=0$ (resp.
$J(R_1,..R_p)(x)=0$) . The zeros of $[\D]$ correspond to the
points of ${\cal{NL}}(\varphi)$, and the zeros of $J$ are the
singular points of $F$. We generalise Daigle's jacobian formula
for $(\CC,+)$-actions ([Da]).

\begin{prop} \label{Daigle}
Let $\varphi$ be an algebraic $\CP$-action on $\CN$ satisfying
condition $(H)$. Then there exists an invariant polynomial $E$
such that $[\D]=E \times J$.
\end{prop}
From a geometric viewpoint, this means that ${\cal{NL}}(\varphi)$
is the union of an invariant hypersurface and of the singular set
of $F$. In particular $E$ is constant if $codim \;
{\cal{NL}}(\varphi) \geq 2$.

\begin{thh} \label{Triv}
Let $\varphi$ be an algebraic $\CP$-action on $\CN$ satisfying condition
$(H)$. If $E$ is constant and $F$ is quasi-fibered, then $\varphi$
is trivial.
\end{thh}
Therefore the assumption "quasi-fibered" correspond to some
regularity in the way that $F$ fibres the orbits. In particular
the action is trivial if $F$ defines a topological quotient, i.e.
if $F$ is smooth surjective and separates the orbits.

\begin{cor} \label{Triv2}
Let $\varphi$ be an algebraic $(\CC,+)$-action on $\CN$ satisfying condition
$(H)$. If $F$ is quasi-fibered, there exists a polynomial $P$ such that
$\varphi$ is conjugate to the action $\varphi'(t;x_1,..,x_n)=(x_1 + tP(x_2,..,x_n),x_2,..,x_n)$.
\end{cor}

\begin{cor} \label{Triv3}
Every algebraic $(\CC ^{n-1},+)$-action $\varphi$ on $\CN$ such
that $codim \; {\cal{NL}}(\varphi)\geq 2$ is trivial. In
particular $\varphi$ is free.
\end{cor}
We end up with counter-examples illustrating the necessity of the
conditions of theorem \ref{Triv} and its corollaries.

\section{Proof of Proposition \ref{torsion}}

In this section, we establish the first proposition announced in
the introduction in two steps. First we describe a TR-exact 1-form
$\omega$ on every generic fibre of $F$. Second we "glue" all these
descriptions by using the uncountability of complex numbers. To
that purpose, we use the following definitions. For any ideal $I$,
we denote by $I\Omega^1(\CN)$ the space of polynomial 1-forms with
coefficients in $I$. We introduce the equivalence relation:
$$\omega \simeq \omega'\;[I] \Longleftrightarrow \omega - \omega '
\in d\Omega^0(\CN)+ \sum \Omega^0(\CN)df_i+I\Omega^1(\CN)$$ This
equivalence is compatible with the structure of $\CC[F]$-module
given by the natural multiplication, since $d\Omega^0(\CN)+ \sum
\Omega^0(\CN)df_i$ and $I\Omega^1(\CN)$ are both $\CC[F]$-modules.

\begin{lem} \label{torsion2}
Let $F^{-1}(y)$ be a non-empty non-critical fibre of $F$, where
$y=(y_1,...,y_q)$. A polynomial 1-form $\omega$ is exact on
$F^{-1}(y)$ if and only if there exists a
polynomial $R$ and some polynomial 1-forms $\eta_1,..,\eta_q$ such that
$\omega = dR +\sum_i(f_i - y_i)\eta_i$.
\end{lem}
\dem Since $\omega$ is exact on $F^{-1}(y)$, it has an holomorphic
integral $R$ on this fibre. Since $F^{-1}(y)$ is a smooth affine
variety, $R$ is a regular map by Grothendieck's Theorem
(\cite{Di}, p. 182). In other words, $R$ is the restriction to
$F^{-1}(y)$ of a polynomial, which will also be denoted by $R$.
The $(q+1)$-form $(\omega - dR) \wedge df_1 \wedge ..\wedge df_q$
vanishes on $F^{-1}(y)$. Since $F^{-1}(y)$ is non-critical, $(f_1
- y_1),.., (f_q - y_q)$ define a local system of parametres at any
point of $F^{-1}(y)$. So the ideal $((f_1 - y_1),...,(f_q - y_q))$
is reduced and we get: $$ (\omega - dR) \wedge df_1 \wedge
..\wedge df_q \equiv 0 \;[f_1-y_1,..,f_q-y_q] $$ The $q$-form
$df_1 \wedge ... \wedge df_q$ never vanishes on $F^{-1}(y)$. By de
Rham Lemma (\cite{Sai}), there exist some polynomials $\alpha_i$
and some polynomial 1-forms $\eta_i$ such that: $$ \omega - dR=
\sum_{i=1} ^q \alpha_idf_i+ \sum_{i=1} ^q (f_i -y_i)\eta_i $$
which can be rewritten as: $$ \omega = d\left(R + \sum_{i=1} ^q
\alpha_i(f_i-y_i) \right) + \sum_{i=1} ^q (f_i
-y_i)(\eta_i-d\alpha_i) $$ \qed {\bf{Proof of Proposition
\ref{torsion}}:} Let $\omega$ be a TR-exact 1-form. Let us show
there exists a non-zero polynomial $P$ such that $P(F)\omega
\simeq 0\;[(0)]$. By lemma \ref{torsion2}, there exists a
non-empty Zariski open set $U$ in $\CQ$ such that, for any
$y=(y_1,..,y_q)$ in $U$: $$ \omega \simeq 0\; [f_1-y_1,..,f_q-y_q]
$$ We proceed to an elimination of $f_1-y_1,..,f_q-y_q$. For any
point $y=(y_{i+1},..,y_q)$ in $\CC^{q-i}$, we denote by $I_i(y)$
the following ideal: $$ I_i(y)=(f_{i+1}-y_{i+1},..,f_q-y_q) $$ By
convention, $\CC^0$ is the space reduced to a point, and
$I_q(y)=(0)$. Let us show by induction on $i\leq q$ the following
property: \\ \ \\ {\em There exists a non-empty Zariski open set
$U_i$ in $\CC ^{q-i}$ such that, for any point $y$ in $U_i$, there
exists a non-zero polynomial $P$ in $\CC[t_1,..,t_i]$ for which
$P(f_1,..,f_i)\omega \simeq 0 \;[I_i(y)]$}. \\ \ \\ This property
is true for $i=0$. Assume it holds to the order $i<q$, and let
$U_i$ be such a Zariski open set. We may assume that $U_i$ is a
principal open set, i.e $U_i= \{f(y)\not=0 \}$. Write
$f=\sum_{k\leq s} f_k(t_{i+2},..,t_q)t_{i+1} ^k$, and set
$U_{i+1}=\{f_s(y')\not=0 \}$. Let $y'=(y_{i+2},...,y_q)$ be a
point in $U_{i+1}$. For any $z$ such that $f(z,y')\not=0$, the
point $y=(z,y')$ belongs to $U_i$. By induction, there exist a
non-zero polynomial $P^z$ and a polynomial 1-form $\eta ^z$ such
that: $$ P^z(f_1,..,f_i)\omega \simeq (f_{i+1}-z)\eta ^z
\;[I_{i+1}(y')] $$ For any such $z$, fix a 1-form $\eta^z$
satisfying this equivalence. The system $\{\eta^z \}$ thus
obtained is an uncountable subset of $\Omega^1(\CN)$. Since
$\Omega^1(\CN)$ has countable dimension, these forms cannot be
linearly independent. There exist some distinct values
$z_1,..,z_m$ and some non-zero constants $(\beta_1,..,\beta_m)$
such that: $$ \beta_1 \eta ^{z_1}+ ... + \beta_m\eta ^{z_m} = 0 $$
Since the equivalence relation is compatible with the structure of
$\CC[F]$-module, we get with the previous relations: $$ \left
(\sum_{j=1} ^m \beta_j P^{z_j}(f_1,..,f_i)\prod _{k\not=j}
(f_{i+1} - z_k) \right ) \omega \simeq 0 \;[I_{i+1}(y')]. $$ None
of the $\beta_j$ (resp. $P^{z_j}$) is zero by construction. Thus
the polynomial $\tilde{P}$: $$ \tilde{P}=\sum_{j=1} ^m \beta_j
P^{z_j}(t_1,..,t_i)\prod _{k\not=j}(t_{i+1} - z_k) $$ is non-zero,
and satisfies the relation $\tilde{P}(f_1,...,f_{i+1})\omega
\equiv 0 \;[I_{i+1}(y')]$. Since we can perform this process for
any point $y'$ in $U_{i+1}$, the induction is proved. \qed

\section{A factorisation lemma}

In this section, we prove an extension of the first Bertini's theorem and
Stein's factorisation theorem ([Sh], p. 139 and [Ha], p. 280) to the case of
reducible varieties. This result is certainly well-known but I could
not find a proper reference for it. So I prefer to give a proof of it,
based on Zariski's Main Theorem.

\begin{lem}
Let $F:X\rightarrow Y$ be a dominating morphism of complex affine varieties,
where $X$ is equidimensional and $Y$ is irreducible. Let $R$ be a regular
map on $X$. Assume that:
\begin{itemize}
\item{The fibres of $F$ are generically connected,}
\item{The restriction of $F$ to any irreducible component of $X$ is dominating,}
\item{The map $G=(F,R)$ is everywhere singular on $X$.}
\end{itemize}
Then $R$ coincides on a dense open set of $X$ with $\alpha(F)$, where $\alpha$ is a rational map on $Y$. In this case, $R$ is said to factor through $F$.
\end{lem}
\dem Since the map $G:X \rightarrow Y \times \CC$ is everywhere singular, $G$
cannot be dominating. So there exists an element $P$ of $\CC[Y][t]$ such that
$P(F,R)=0$ on $X$. Note that $P$ has degree $>0$ with respect to $t$, because
$F$ is a dominating map. Under the previous assumptions, there exists
a Zariski open set $U$ in $Y$ such that:
\begin{itemize}
\item{For any irreducible component $X'$ of $X$, $U$ is contained in $F(X')$,}
\item{For any point $y$ in $U$, $F^{-1}(y)$ is connected,}
\item{For any point $y$ in $U$, the polynomial $P(y,t)$ is non-zero.}
\end{itemize}
Let $y$ be a point in $U$. Since $P(y,R)=0$ on $F^{-1}(y)$, $R$ is
locally constant on $F^{-1}(y)$. Since $R$ is regular and
$F^{-1}(y)$ is connected, $R$ is constant on $F^{-1}(y)$. So we
can define the correspondence $\alpha: U\rightarrow \mathbb C$
that maps any point $y$ of $U$ to the unique value that takes $R$
on $F^{-1}(y)$. Consider its graph: $$ Z=\{(y,\alpha(y)),y \in U
\} $$ If $X'$ is an irreducible component of $X$, then $Z$
coincides with $G(X'\cap F^{-1}(U))$. So $Z$ is constructible for
the Zariski topology, and $\overline{Z}$ is irreducible. Therefore
$\overline{Z}$ defines in $Y\times\mathbb C$ a rational
correspondence from $Y$ to $\mathbb C$ in the sense of Zariski
(\cite{Mu1}, pp. 29-51). By Zariski's Main Theorem, $\alpha$
coincides with a rational map on $Y$. Let $U'$ be an open set
contained in $U$ where $\alpha$ is regular. Then $F^{-1}(U')$ is a
dense open subset of $X$. Moreover $R$ and $\alpha (F)$ coincide
on $F^{-1}(U')$ by construction. \qed

\section{Blowing-downs and primitive mappings}

In this section, we give some properties of blowing-downs and
primitive mappings. For this class of maps, we will establish a
{\em division lemma} (see section 5) that is the key-point for the
proof of theorem \ref{quasi}. Let $F$ be a polynomial dominating
map from $\CN$ to $\CQ$, and let $S(F)$ be its set of singular
points. We introduce the following sets: $$
\begin{array}{ccl}
B(F)& = & \{y \in \mathbb C ^q,\mbox{\it $F^{-1}(y)$ is non-empty and not
connected}\} \ \\
E(F) & = & \mbox{\it Union of blowing-downs of $F$} \ \\
I(F) & = &\{y \in \mathbb C^q,\mbox{\it $F^{-1}(y)$ is empty}\}
\end{array}
$$ Let $H$ be the GCD of all $q$-minors of $dF$, and set: $$
\omega _F = \frac{df_1\wedge .. \wedge df_q}{H} $$ Note that for
all polynomials $P$ and $R$, we have $P(F) dR \wedge \omega _F =
d(P(F) R) \wedge \omega _F$. Since the sets $B(F),E(F),I(F)$ are
all constructible for the Zariski topology, it makes sense to
consider their codimensions. Recall that $F$ is primitive if its
fibres are 0-generically connected and 1-generically non-empty,
i.e. $codim\; B(F) \geq 1$ and $codim\; I(F)\geq 2$.

\begin{prop} \label{prim}
A polynomial map $F:\CN \rightarrow \CQ$ is primitive if
and only if any polynomial $R$ such that $dR \wedge \omega_F = 0$
belongs to $\CC[F]$.
\end{prop}
\dem Assume that $F$ is primitive. Let $R$ be a polynomial such that
$dR \wedge \OF =0$. Then the map $G=(F,R)$ is everywhere singular.
Since the generic fibres of $F$ are connected, $R$ factors through
$F$ by the factorisation lemma. Let us set:
$$
R=b(F)/a(F)
$$
where $a,b$ are
relatively prime. Let us show by absurd
that $a$ is constant. Assume not, and let $a'$ be an irreducible factor of
$a$. For any point $y$ in $V(a') - I(F)$, there exists a point $x$ such that
$F(x)=y$, which implies that $a(y)R(x)=b(y)=0$. So $b$ vanishes on
$V(a') - I(F)$. Since $I(F)$ has
codimension $\geq 2$ in $\CN$, $V(a') - I(F)$ is dense in $V(a')$ and
$b$ vanishes on $V(a')$. By Hilbert's
Nullstellensatz, $a'$ divides $b$, contradicting the fact that $a$
and $b$ are relatively prime.
Thus $a$ is constant and $R$ belongs to $\CC[F]$.

Assume now that any polynomial $R$ such that $dR\wedge \OF=0$ belongs to
$\CC[F]$. The $q$-form $\OF$ is obviously non-zero, and the polynomials
$f_i$ are algebraically independent. So $F$ is a dominating map. Let us
prove first that $codim(B(F))\geq 1$. By Bertini first
theorem ([Sh], p. 139), it suffices to show that $\CC(F)$
is algebraically closed in $\CC(x_1,...,x_n)$. Let $R$ be
a rational fraction that is algebraic over $\CC(F)$. Let
$P(z,t_1,..,t_q)=\sum_{k\leq s} a_k(t_1,..,t_q)z^k$ be a
nonzero polynomial such that $P(R,f_1,..,f_q)=0$. We choose
$P$ of minimal degree with respect to $z$. Since
$P(R,f_1,..,f_q)=0$, the denominator of $R$ divides
$a_s(F)$. By derivation and wedge product, we get:
$$
\frac{\partial  P}{\partial z}(R,f_1,..,f_q)dR \wedge \OF =0
$$
Since $P$ has minimal degree, $dR \wedge \OF=0$ and
$d(a_s(F)R) \wedge \OF=0$. As $a_s(F)R$ is a polynomial,
it belongs to $\CC[F]$ and $R$ lies in $\CC(F)$.

Let us show by absurd that $codim(I(F))\geq 2$. Assume not, and
let $C=V(f)$ be a codimension 1 irreducible component of
$\overline{I(F)}$, where $f$ is reduced. Since the intersection
$V(f) \cap F(\CN)$ has codimension $\geq 2$, there exists a
polynomial $P$ vanishing on $V(f) \cap F(\CN)$ and not divisible
by $f$. The function $P(F)$ vanishes on $V(f(F))$. By Hilbert's
Nullstellensatz, there exists an integer $n$ such that $P^n(F)$ is
divisible by $f(F)$. The function $P^n/f$ is rational
non-polynomial, and $R=P^n(F)/f(F)$ belongs to $\CX$. Since $R$
satisfies the equation $dR \wedge \OF =0$, $R$ belongs to $\mathbb
C[F]$, hence a contradiction. \qed For $q=1$, a mapping $F$ is
primitive if and only if its generic fibres are connected. Indeed
any non-constant polynomial map from $\CN$ to $\CC$ has to be
surjective. In this way, the definition of primitive mapping
extends the notion of primitive polynomial ([D-P]). \\ \ \\
{\bf{Exemple 1:}} The polynomial $F(x,y)=x^2$ is not primitive
because its generic fibres are not connected. Note that $dx\wedge
d(x^2)=0$, but $x$ does not belong to $\CC[x^2]$. \\ \ \\
{\bf{Exemple 2:}} Consider the mapping $F: \mathbb C^3 \rightarrow
\mathbb C^2, (x,y,z)\mapsto (x,xy)$. The function $y$ satisfies
the relation $dy\wedge dx \wedge d(xy) =0$ but does not belong to
$\mathbb C[x,xy]$. So $F$ is not a primitive mapping although its
generic fibres are connected. The obstruction lies in the fact
that $\overline{I(F)}=\{(y_1,y_2), y_1=0\}$, so $codim(I(F))=1$.
\\ \ \\ {\bf{Exemple 3:}} Consider the mapping $F: \mathbb C^3
\rightarrow \mathbb C^2, (x,y,z)\mapsto (xy,zy)$. It is easy to
see that $F$ is onto and that its generic fibres are isomorphic to
$\mathbb C^*$. So $F$ is a primitive mapping.
\\ \
\\ Recall that a blowing-down is an hypersurface of $\CN$ that is
mapped by $F$ to a set of codimension $\geq 2$. For instance, the
plane $\{y=0\}$ in $\CC^3$ is a blowing-down of the map
$F(x,y,z)=(xy,zy)$.

\begin{prop} \label{Blow}
Any blowing-down of $F$ is contained in $S(F)$.
\end{prop}
\dem Let $V$ be a blowing-down of $F$, and let $W$ denote
the Zariski closure of $F(V)$. Then $W$ is irreducible and
there exists a dense open set $W'$ of $W$, consisting only
of smooth points of $W$ and containing $F(V)$. So $V'=
F^{-1}(W')\cap V$ is a dense open set of $V$.
For any smooth point $x$ in $V'$, the differential of
the restriction of $F$ to $V$ has rank
$\leq dim W' \leq q-2$. The differential $dF(x)$
maps the hyperplane $T_x V$ to a space of dimension
$\leq q-2$.
So $dF(x)$ maps $\CN$ to a space of dimension $\leq q-1$,
and $F$ is singular at $x$. Since any smooth point of $V'$
is a singularity of $F$ and $S(F)$ is closed, we have the
inclusion $V\subset S(F)$.
\qed

\section{The division lemma}

In this section, we are going to establish the essential tool for the proof
of theorem \ref{quasi}. Let $\omega$ be a TR-exact 1-form $\omega$.
By proposition \ref{torsion}, there exists a non-zero polynomial $P$ in $\CT$,
and some polynomials $R,a_1,...,a_q$ in $\CX$ such that:
$$
P(F)\omega = dR+ a_1df_1+...+a_qdf_q
$$
By using the wedge product with $\omega_F$, we get :
$$
dR \wedge \omega_F=P(F)\omega \wedge \omega_F\equiv 0\;[P(F)]
$$
Assume there exist some polynomials $S,b_1,..,b_q$ such that
$\omega=dS + \sum_i b_idf_i$. By an obvious computation, we
get $\omega \wedge \OF = dS \wedge \OF$ and
$d(R-P(F)S)\wedge \OF=0$. Since $F$ is primitive, there
exists a polynomial $A$ such that $R=A(F) + P(F)S$.

More generally, let $R$ be a polynomial satisfying the equation
$dR \wedge \omega_F\equiv 0\;[P(F)]$. $R$ is said to be {\em
${\cal{E}}$-divisible} by $P(F)$ if there exist some polynomials $A$ and
$S$ such that $R=A(F) + P(F)S$. In this section we are going
to determine under which conditions a polynomial $R$ satisfying
this equation is ${\cal{E}}$-divisible by $P(F)$. \\ \ \\
{\bf{Division lemma}}
{\it Let $F$ be a primitive mapping from $\CN$ to $\CQ$. Let
$P$ be an element of $\CT$, and $R$ a polynomial in $\CX$
satisfying
the equation $dR \wedge \omega_F \equiv 0\;[P(F)]$. Assume
that:
\begin{itemize}
\item{$V(P)\cap B(F)$ has codimension $\geq 2$ in $\CQ$,}
\item{$V(P(F)) \cap E(F)$ has codimension $\geq 2$ in $\CN$,}
\item{$V(P)\cap I(F)$ has codimension $\geq 3$ in $\CQ$.}
\end{itemize}
Then $R$ is ${\cal{E}}$-divisible by $P(F)$.}

\subsection{The weak division lemma}

In this subsection, we are going to establish a weak version of the
division lemma. A polynomial $R$ is said to be {\em weakly ${\cal{E}}$-divisible}
by $P(F)$ if there exists a polynomial $B$ coprime to $P$ such
that $B(F)R$ is ${\cal{E}}$-divisible by $P(F)$. \\ \ \\
{\bf{Weak division lemma}}
{\it Let $F$ be a primitive mapping from $\CN$ to $\CQ$. Let
$P$ be an irreducible polynomial of $\CT$. Let $R$ be a polynomial
in $\CX$
satisfying
the equation $dR \wedge \omega_F \equiv 0\;[P(F)]$. Assume
that:
\begin{itemize}
\item{$V(P)\cap B(F)$ has codimension $\geq 2$ in $\CQ$,}
\item{$V(P(F)) \cap E(F)$ has codimension $\geq 2$ in $\CN$,}
\end{itemize}
Then $R$ is weakly ${\cal{E}}$-divisible by $P(F)$.} \\ \ \\
The proof splits in two steps. Consider a polynomial
$R$ satisfying the equation
$dR \wedge \omega_F \equiv 0\;[P(F)]$. First we show that its
restriction to $V(P(F))$ factors through $F$. So there exist
two polynomials $A,B$, with $B$ coprime to $P$, such that
$B(F)R-A(F)$ vanishes on $V(P(F))$. If $h_1 ^{n_1}.. h_r ^{n_r}$
is the irreducible decomposition of $P(F)$ in $\CX$, then
$h_1 .. h_r $ divides $B(F)R-A(F)$. Second we prove that
every factor $h_i$ divides $B(F)R-A(F)$ with multiplicity $\geq n_i$.

\begin{lem}
Let $P$ be an irreducible polynomial in $\CT$. Let $h$ be an
irreducible factor of $P(F)$. Let $R$ be a polynomial
satisfying the equation $dR \wedge \OF \equiv 0\;[h]$.
Then the map $G:V(h)\rightarrow V(P)\times \CC, \;x\mapsto (F(x),R(x))$
is everywhere singular.
\end{lem}
\dem It suffices to show that the collection of 1-forms
$dR,dh,df_1,..,df_q$ has rank $\leq q$ at any point $x$ of $V(h)$.
We are going to check that whenever you choose $q+1$ forms in this
collection, their wedge product is divisible by $h$. Consider the
first case, when this wedge product contains all the forms
$df_1,..,df_q$. Then it is either equal to $dR \wedge df_1 \wedge
..\wedge df_q$ or to $dh \wedge df_1 \wedge ..\wedge df_q$. By
assumption $dR \wedge df_1 \wedge ..\wedge df_q$ is divisible by
$h$. To see that the second one is divisible by $h$, factor
$P(F)=Qh^m$, where $Q$ is coprime to $h$ and $m\geq 1$. By wedge
product, we get: $$ d P(F) \wedge df_1 \wedge ..\wedge df_q =
mh^{m-1}Q dh \wedge df_1 \wedge ..\wedge df_q + h^{m}dQ  \wedge
df_1 \wedge ..\wedge df_q =0 $$ This yields $Q dh \wedge df_1
\wedge ..\wedge df_q\equiv 0[h]$. Since $Q$ is coprime to $h$, we
find: $$ dh \wedge df_1 \wedge ..\wedge df_q\equiv 0[h] $$
Consider now the second case, when $dR$ and $dh$ appear in the
wedge product. Assume first that $q>1$. Up to a reordering of the
forms $df_i$, we may assume that this wedge product is equal to
$dR \wedge dh \wedge df_2 \wedge ..\wedge df_q$. Since $P(F)=
Qh^m$ where $Q$ is coprime to $h$, we get by derivation: $$
d\{P(F)\}=\sum_{i=1} ^q \frac{\der P}{\der t_i}(F)df_i \equiv 0 \;
[h^{m-1}] $$ By wedge product, we find: $$ \frac{\der P}{\der
t_1}(F)H\OF = \frac{\der P}{\der t_1}(F)df_1 \wedge ..\wedge df_q
= d\{P(F)\}\wedge df_2 \wedge ..\wedge df_q \equiv 0 \; [h^{m-1}]
$$ By construction, the coefficients of $\OF$ have no common
factors. Thus $h^{m-1}$ divides $\der P /\der t_1 (F)H$. Then
write: $$ \frac{\der P}{\der t_1}(F)H dR \wedge \OF= dR \wedge
d\{P(F)\}\wedge df_2 \wedge ..\wedge df_q = dR \wedge
d\{Qh^m\}\wedge df_2 \wedge ..\wedge df_q $$ Since $dR \wedge \OF$
is divisible by $h$, we get: $$ dR \wedge d\{Qh^m\}\wedge df_2
\wedge ..\wedge df_q\equiv 0\; [h^m] $$ which leads to: $$
mQh^{m-1}dR \wedge dh \wedge df_2 \wedge ..\wedge df_q\equiv 0\;
[h^m] $$ Since $Q$ is coprime to $h$, we deduce: $$ dR \wedge dh
\wedge df_2 \wedge ..\wedge df_q\equiv 0\; [h] $$ If $q=1$, we do
the same computation and forget the wedge product with $df_2
\wedge ..\wedge df_q$. \qed

\begin{lem} \label{divv}
Let $P$ be an irreducible polynomial in $\CT$. Let
$h_1 ^{n_1}.. h_r ^{n_r}$ be the irreducible
decomposition of $P(F)$ in $\CX$. Let $R$ be
a polynomial such that $dR \wedge \OF
\equiv 0\; [h_1.. h_r]$. Assume that:
\begin{itemize}
\item{$V(P)\cap B(F)$ has codimension $\geq 2$ in $\CQ$,}
\item{$V(P(F)) \cap E(F)$ has codimension $\geq 2$ in $\CN$.}
\end{itemize}
Then there exist two polynomials $A,B$, where $B$ is coprime to
$P$, such that $B(F)R -A(F)$ is divisible by $h_1.. h_r$.
\end{lem}
\dem By the previous lemma applied to all the irreducible
components of $V(P(F))$, we can see that the map: $$ G:
V(P(F))\rightarrow V(P)\times \CC , \; x\mapsto (F(x),R(x)) $$ is
singular. Since $V(P(F))\cap E(F)$ has codimension $\geq 2$, none
of the hypersurfaces $V(h_i)$ is a blowing-down. So $F$ maps every
$V(h_i)$ densely on $V(P)$. Since $V(P)\cap B(F)$ has codimension
$\geq 2$, the generic fibres of $F:V(P(F))\rightarrow V(P)$ are
connected. By the factorisation lemma, there exists a rational map
$\alpha$ on $V(P)$ such that $R=\alpha(F)$ on $V(P(F))$. Write
$\alpha$ as $A/B$, where $B$ is coprime to $P$. The polynomial
$B(F)R -A(F)$ vanishes on $V(P(F))$. By Hilbert's Nullstellensatz,
it is divisible by $h_1.. h_r$. \qed {\bf{Proof of the weak
division lemma:}} Let $P$ be an irreducible polynomial in $\CT$.
Let $h_1 ^{n_1}.. h_r ^{n_r}$ be the irreducible decomposition of
$P(F)$ in $\CX$. Let $R$ be a polynomial such that $dR \wedge \OF
\equiv 0\; [P(F)]$. Then $R$ satisfies the equation: $$ dR \wedge
\OF \equiv 0\; [h_1 .. h_r ] $$ By the previous lemma, there exist
some polynomials $A,B$, where $B$ is coprime to $P$, such that
$S=B(F)R-A(F)$ is divisible by $h_1..h_r $. Factor $S$ as $S_0h_1
^{k_1}.. h_r ^{k_r}$, where $S_0$ is coprime to each $h_i$. Let us
show by absurd that $k_i\geq n_i$ for any $i$.

Assume there exists an index $i$ such that $k_i/n_i <1$. Let $i_0$
be an index for which the ratio $k_i/n_i$ is minimal, and let
$u/v$ be its irreducible decomposition. By construction, we have
$0<u/v<1$. The function: $$ L=S^v /P(F)^u=S_0 ^v h_1 ^{vk_1 -
un_1}.. h_r ^{vk_r - un_r} $$ is polynomial, since $u/v \leq
k_i/n_i \Rightarrow vk_i - un_i\geq 0$. Moreover $L$ satisfies the
equation $dL \wedge \OF \equiv 0\;[h_1..h_r]$. Indeed if $vk_i -
un_i>0$, then $L$ is divisible by $h_i$ and $L=L_i h_i$. We set
$P(F)=P_ih_i ^{n_i}$, where $P_i$ is coprime to $h_i$. By an easy
computation, we get: $$ dP(F) \wedge \OF = P_i n_i h_i ^{n_i
-1}dh_i \wedge \OF + h_i ^{n_i}dP_i \wedge \OF =0 $$ Since $P_i$
is coprime to $h_i$, we deduce $dh_i \wedge \OF \equiv 0\;[h_i]$,
and this implies: $$ dL \wedge \OF =L_idh_i \wedge \OF + h_idL_i
\wedge \OF \equiv 0\;[h_i] $$ If $vk_i - un_i=0$, set $S=S_i h_i
^{k_i}$. By derivation and wedge product, we get: $$ SdL \wedge
\OF= S_i h_i ^{k_i}dL \wedge \OF = vLdS \wedge \OF $$ By an easy
computation, we obtain: $$ dS \wedge \OF =B(F)dR \wedge \OF \equiv
0 \;[h_i ^{n_i}] $$ which implies: $$ S_i dL \wedge \OF \equiv 0
\;[h_i ^{n_i -k_i}] $$ Since $n_i -k_i >0$ and $S_i$ is coprime to
$h_i$, we deduce $dL \wedge \OF \equiv 0 \;[h_i]$. Thus $dL \wedge
\OF$ is divisible by $h_1..h_r$. By lemma \ref{divv}, there exist
two polynomials $A',B'$, where $B'$ is coprime to $P$, such that
$B'(F)L - A'(F)\equiv 0\;[h_1..h_r]$.

Let us show by absurd that $vk_i - un_i=0$ for any $i$. Assume
that $h_i$ divides $L$. By the previous relation, $h_i$
divides $A'(F)$. Since $V(h_i)$ is not a blowing-down
and $P$ is irreducible, $A'$ is divisible by $P$, which
implies:
$$
B'(F)L\equiv 0\;[h_1..h_r]
$$
Since none of the $V(h_j)$ are blowing-downs and every $h_j$
divides $P(F)$, every $h_j$ is coprime to $B'(F)$. So $L$
is divisible by $h_1..h_r$, contradicting its construction.

Since $vk_i - un_i=0$, $v$ divides $n_i$ for any $i$. As
$0<u/v<1$, $v$ is strictly greater than 1 and $P(F)=T^v$,
where $T$ belongs to $\CX$. This implies:
$$
d\{P(F)\} \wedge \OF = vT^{v-1} dT \wedge \OF =0
$$
Since $F$ is primitive, $T$ belongs to $\CC[F]$ by proposition \ref{prim}.
Therefore $P$ is the $v^{th}$ power of some polynomial, which contradicts
the irreducibility of $P$.
\qed

\subsection{Proof of the division lemma}

Let $R$ be a polynomial satisfying the equation $dR \wedge \OF
\equiv 0\; [P(F)]$. From an analytic viewpoint, the weak division
lemma asserts that $R$ coincides on $V(P(F))$ with $\alpha(F)$,
where $\alpha$ is a rational function on $V(P)$. In order to prove
the division lemma, we are going to show that $\alpha$ is regular
if $V(P)\cap I(F)$ has codimension $ \geq 3$. In other words, we
are going to eliminate the "poles" of $\alpha$.

Recall that an ideal $I$ in a local ring $R$ is ${\cal{M}}$-primary if
$I$ contains some power of the maximal ideal ${\cal{M}}$ of $R$. We
denote by ${\cal{O}}_{\CQ,y}$ the ring of germs of regular functions
at the point $y$ in $\CQ$. For simplicity, we set:
$$
\CC[[X]]=\CC[[x_1,..,x_n]] \quad \mbox{and} \quad \CC[[T]]=\CC[[t_1,..,t_q]]
$$
\begin{lem} \label{div2}
Let $I=(g_1,..,g_n)$ be an ${\cal{M}}$-primary ideal in
$\CC[[X]]$. If the
classes of the formal series $\{e_1,..,e_{\mu}\}$ form
a basis of the
vector space $\CC[[X]]/I$, then $\{e_1,..,e_{\mu}\}$ is
a basis of
the $\CC[[g_1,..,g_n]]$-module $\CC[[X]]$.
\end{lem}
\dem Since $(g_1,..,g_n)$ is ${\cal{M}}$-primary, $\CC[[X]]$ is a finitely
generated $\CC[[g_1,..,g_n]]$-module ([Ab]). By Nakayama lemma
([Sh], p. 283), $\{e_1,...,e_{\mu}\}$ forms a minimal set
of generators of this module. Let us show by absurd that
$e_1,...,e_{\mu}$ are $\CC[[g_1,..,g_n]]$-linearly
independent.

Assume there exist some formal series $a_i(y_1,..,y_n)$, not all
equal to zero, such that $\sum_k a_k(g_1,..,g_n)e_k=0$. Up to a
linear change of
coordinates on $y_1,..,y_n$, which is equivalent to replacing
$g_1,..,g_n$ by another set of formal series generating the
same ideal, we may assume there exists an index $i$ for
which $a_i(y_1,0,..,0)\not=0$. By setting $a_i(x_1,0,..,0)
=b_i(x_1)$, we find:
$$
b_1(g_1)e_1 + .. + b_{\mu}(g_1)e_{\mu}\equiv 0 \;[g_2,..,g_n]
$$
Let $m$ be the minimum of the orders of all formal series
$b_1,...,b_{\mu}$. Then $b_i(x_1) =x_1 ^m c_i(x_1)$ for any
$i$, and $c_i(0)\not=0$ for at least one of them. Thus we get:
$$
g_1 ^m \{ c_1(g_1)e_1 + .. + c_{\mu}(g_1)e_{\mu}
\}\equiv 0 [g_2,..,g_n]
$$
Since $(g_1,..,g_n)$ is ${\cal{M}}$-primary, $g_1,..,g_n$ is a
regular sequence ([Sh], p. 227) and $g_1$ is not a zero-divisor
modulo $[g_2,..,g_n]$. We deduce:
$$
c_1(0)e_1 + .. + c_{\mu}(0)e_{\mu}\equiv 0 [g_1,g_2,..,g_n]
$$
So $c_1(0)=..=c_{\mu}(0)=0$, hence contradicting the fact that not all
$c_i(0)$ are zero.
\qed

\begin{lem} \label{div3}
Let $y$ be a point in $\CQ$ such that the fibre $F^{-1}(y)$
is non-empty
of dimension $(n-q)$. Let $P,B,A$ be three elements of $\CT$ such
that $A(F)$ belongs to the ideal $(P(F),B(F))\CX$.
Then $A$ belongs to $(P,B){\cal{O}}_{\CQ,y}$.
\end{lem}
\dem Let $x$ be a point in $F^{-1}(y)$ where the fibre has local
dimension $(n-q)$. For simplicity, we may assume $x=0$ and $y=0$.
There exists a $q$-dimensional vector space, defined by some
linear equations $l_1,...,l_{n-q}$ and intersecting locally
$F^{-1}(0)$ only at 0. By Ruckert's Nullstellensatz ([Ab]), the
ideal $(f_1,..,f_q,l_1,..,l_{n-q})$ is ${\cal{M}}$-primary in the
ring $\CC[[X]]$. Let $\{e_1,..,e_{\mu}\}$ be a basis of the vector
space $\CC[[X]]/(f_1,..,f_q,l_1,..,l_{n-q})$ such that $e_1=1$. By
lemma \ref{div2}, $\{e_1,..,e_{\mu}\}$ is a basis of the
$\CC[[f_1,..,l_{n-q}]]$-module $\mathbb C[[X]]$. Let $R,S$ be two
polynomials in $\CX$ such that $A(F)=P(F)R+Q(F)S$. If
$R_1(f_1,..,l_{n-q})$ and $S_1(f_1,..,l_{n-q})$ denote their first
coordinate in the basis $\{e_1,..,e_{\mu}\}$, we get: $$
P(F)R_1(f_1,..,l_{n-q}) + B(F)S_1(f_1,..,l_{n-q})=A(F) $$ After
reduction modulo $l_1,..,l_{n-q}$, this implies: $$ P(F)R_1(F,0) +
B(F)S_1(F,0)=A(F) $$ Thus $A$ belongs to the ideal $(P,B)\mathbb
C[[T]]$. Since ${\cal{O}}_{\CQ,0}$ is a Zariski ring and $\mathbb
C[[T]]$ is its ${\cal{M}}$-adic completion, we get $(P,B)\CC[[T]]
\cap {\cal{O}}_{\CQ,0}=(P,B){\cal{O}}_{\CQ,0}$ (\cite{Ma1}, pp.
171-172). So $A$ belongs to $(P,B){\cal{O}}_{\CQ,0}$. \qed

\begin{lem} \label{div4}
Let $P,B,A$ be three polynomials in $\CT$ such that $A(F)$
belongs to
$(P(F),B(F))\CX$. If $V(P(F),B(F))$ has codimension $\geq 2$
and $V(P(F))\cap I(F)$ has codimension $\geq 3$, then $A$
belongs to
$(P,B)\CT$.
\end{lem}
\dem This lemma is obvious if $V(P,B)$ is empty. We assume it is not,
and consider the varieties $X=V(P(F),B(F))$ and $Y=V(P,B)$. By assumption,
$P(F)$ and $B(F)$ are coprime and $X$ is equidimensionnal of codimension
2 in $\CN$. Moreover $P,B$ are coprime and $Y$ is equidimensionnal
of codimension 2 in $\CQ$. As $V(P)\cap I(F)$ has codimension $\geq 3$,
the restriction:
$$
F_R : X\longrightarrow Y, \; x \longmapsto F(x)
$$
is a dominating map. We construct a dense open set $U$ in $Y$ such
that $F^{-1}(y)$ has dimension $(n-q)$ for any $y$ in $U$. Let $X_i$
be any irreducible component of $X$. If $F(X_i)$ has codimension $\geq 3$,
fix a dense open set $U_i$ in $Y$ that does not meet $F(X_i)$.
If $F(X_i)$
has codimension 2, we apply the theorem on the dimension of fibres
to $F_R: X_i \rightarrow \overline{F(X_i)}$. There exists an open
set $V_i$ contained in $F(X_i)$ such that $F^{-1}(y)\cap X_i$ has
dimension $(n-q)$ for any $y$ in $V_i$. If $U'$ is the intersection
of
all $U_i$ and $V'$ is the union of all $V_i$, then $U=U' \cap V'$
is a dense open set in $Y$, and $F^{-1}(y)$ has dimension $(n-q)$
for any $y$ in $U$.

By lemma \ref{div3}, $A$ belongs to $(P,B){\cal{O}}_{\CQ,y}$ for
any $y$ in $U$. This means there exists a polynomial $\beta_y$
such that $\beta_y(y)\not=0$ and $\beta_y A$ belongs to $(P,Q)\CT$.
The zero set of $P,B$ and the $\beta_y$, when $y$ runs through $U$,
has codimension $\geq 3$ since it is contained in $Y-U$. The
ideal $J$ generated by $P,B$ and the $\beta_y$ has depth $\geq 3$.
Since $\CT$ is catenary, $J$ contains a polynomial $\beta$ such
that $P,B,\beta$ is a regular sequence. By construction $\beta A
\equiv 0\;[P,B]$. As $\beta$ is not a zero divisor modulo $(P,B)$,
$A$ belongs to $(P,B)\CT$.
\qed
{\bf{Proof of the division lemma:}} Let $R$ be a polynomial satisfying
the equation $dR \wedge \OF \equiv 0 \; [P(F)]$. Assume that
$V(P)\cap B(F)$ has codimension $\geq 2$, $V(P(F))\cap E(F)$
has codimension $\geq 2$ and $V(P)\cap I(F)$ has codimension
$\geq 3$. By the weak division lemma, there exist two polynomials
$A,B$, where $B$ is coprime to $P$, and a polynomial $S$ such that:
$$
B(F)R - A(F)=P(F)S
$$
Let us show by absurd that $X=V(P(F),B(F))$ has codimension $\geq 2$.
Assume that $X$ contains an hypersurface $V$. Then $F$ maps $V$ to
$Y=V(P,B)$, which codimension is $\geq 2$ since $P$ and $B$ are
coprime. So $V$ is a blowing-down, and this contradicts the
assumption on $V(P(F))\cap E(F)$.

Since $A(F)$ belongs to $(P(F),B(F))\CX$ and $V(P)\cap I(F)$
has codimension $\geq 3$, $A$ belongs to $(P,B)\CT$ by lemma
\ref{div4}. There exist some polynomials $P_1,B_1$ such that
$A=PP_1 + BB_1$. Thus we deduce:
$$
B(F)\{R-B_1(F)\}=P(F)\{S - P_1(F)\}
$$
Since $X=V(P(F),B(F))$ has codimension 2, $P(F)$ and $B(F)$ are
coprime. So $P(F)$ divides $R-B_1(F)$ and the division lemma
is proved.
\qed

\subsection{Proof of theorem \ref{quasi3}}

Let $F$ be a primitive mapping that is
either quasi-fibered or weakly quasi-fibered. By definition, the
following conditions hold:
\begin{itemize}
\item{$B(F)$ has codimension $\geq 2$ in $\CQ$,}
\item{$E(F)$ is empty,}
\item{$I(F)$ has codimension $\geq 3$ in $\CQ$.}
\end{itemize}
Let $\omega$ be a TR-exact 1-form. By proposition \ref{torsion},
there exists a non-zero polynomial $P$, and some polynomials
$R,a_1,..,a_q$ such that: $$ P(F)\omega= dR+ a_1df_1+...+a_qdf_q
$$ By wedge product with $\OF$, we can see that $R$ satisfies the
equation $dR \wedge \OF \equiv 0 \;[P(F)]$. According to the
conditions given above, $V(P)\cap B(F)$ has codimension $\geq 2$
in $\CQ$, $V(P(F))\cap E(F)$ is empty and $V(P)\cap I(F)$ has
codimension $\geq 3$ in $\CQ$. By the division lemma, there exist
some polynomials $A$ and $S$ such that $R=A(F) + P(F)S$. Therefore
a simple calculation yields: $$ P(F)\omega= P(F)dS + \sum_{k=1} ^q
\left ( a_k + S \frac{\partial P}{\partial t_k}(F)+ \frac{\partial
A}{\partial t_k}(F) \right )df_k $$ Let $c_k$ denote the
coefficient of $df_k$ in this sum. Then $\sum_k c_kdf_k$ is
divisible by $P(F)$. If $\omega_0$ is that quotient, we can see:
$$ \omega_0\wedge df_1 \wedge ..\wedge df_q=0 $$ which implies the
second part of the theorem. If now $F$ is quasi-fibered, then it
is non-singular in codimension 1. By De Rham Lemma (\cite{Sai}),
$\omega_0$ can be written as $\sum_k d_k df_k$, where all $d_k$
are polynomials. Therefore $\omega$ is AR-exact. \qed

\section{Recalls on $\CP$-actions}

An algebraic $\CP$-action $\varphi$ on an affine variety $X$
consists of a regular map $\varphi: \CC ^p \times X \rightarrow X$
such that: $$ \forall (u,v) \in \CC^p \times \CC^p, \quad \forall
x \in X, \quad \varphi (u,\varphi (v,x))= \varphi (u+v,x) $$ We
denote by $\CC[X]^{\varphi}$ its ring of invariants, i.e. the
space of regular functions $f$ such that $f\circ \varphi = f$. The
action $\varphi$ can be defined as the composition of $p$ pairwise
commuting algebraic $(\CC,+)$-actions $\varphi_i$. These latter
are the restriction of $\varphi$ to the $i^{th}$ coordinate of
$\CC^p$. To each $\varphi_i$ corresponds the derivation
$\der_i=\varphi_i ^*(d/dt_i)_{t_i=0}$, which enjoys the remarkable
property of being locally nilpotent (see the introduction).
Moreover these derivations commute pairwise. Conversely if
$\{\der_1,...,\der_p\}$ is a system of locally nilpotent
pairwise
commuting derivations, the exponential map: $$
exp(t_1\der_1+..+t_p\der_p)(f)= \sum_{k\geq 0}
\frac{(t_1\der_1+..+t_p\der_p)^k(f)}{k!} $$ defines a morphism of
algebras from $\CC[X]$ to $\CC[X]\otimes \CC[t_1,..,t_p]$. This
morphism induces a regular map $\varphi:\CC ^p \times X
\rightarrow X$ that is an $\CP$-action on $X$. In this case,
$\varphi$ is said to be {\em generated} by $\{\der_1,..,\der_p\}$.

\begin{df}
A commutative $p$-distribution $\D$ is a system of locally nilpotent
pairwise commuting derivations $\der_1,...,\der_p$. Its ring of invariants
$\CC[X] ^{\D}$ is the intersection of the kernels of the $\der_i$ on $\CC[X]$.
\end{df}
If $\varphi$ is generated by $\D$, then $\CC[X] ^{\D}$ is the ring of
invariants of $\varphi$. Indeed, by definition of $\varphi$ via
the exponential map,
a regular function $f$ is invariant by $\varphi$ if and only if
$\der_i(f)=0$ for any $i$.
Recall that the action $\varphi$ is free at $x$ if the stabilizer
of $x$ is reduced to zero, or in other words if the orbit of $x$
has dimension $p$. Let $[\D]$ be the operator defined at the
introduction. We introduce its evaluation at $x$:
$$
[\D](x): (R_1,..,R_p) \longmapsto det((\partial_i(R_j)))(x)
$$
\begin{lem}
Let $\varphi$ be an algebraic $\CP$-action on $X$, and let $\D$ be its
commutative $p$-distribution. Then $\varphi$ is not free at $x$
if and only if $[\D](x)$ is the null map.
\end{lem}
\dem Assume first that $\varphi$ is not free at $x$. Let
$(u_1,..,u_p)$ be a non-zero element of the stabilizer
of $x$. Let $\varphi ^u$ be the $(\mathbb C , +)$-action
defined by $\varphi ^u_t(y)=\varphi _{tu_1,..,tu_p}(y)$.
Starting from the relation $\varphi ^u _1(x)=x$, we get by
an obvious induction that $\varphi ^u _m(x)=x$ for any
integer $m>0$. So $\varphi ^u _t(x)=x$ for any $t$ in $\CC$,
and $x$ is a fixed point of $\varphi ^u$. For any regular function
$R$, we get by derivation:
$$
\sum u_i \partial_i(R)(x)=0
$$
which implies for any $p$-uple $(R_1,..,R_p)$:
$$
[\D](x)(R_1,..,R_p)=det((\partial_i(R_j)))(x)=0
$$
Assume now that $[\D](x)$ is the null map. Let $(\der_i)_x$
be the evaluation map of $\der_i$ at $x$, i.e. the map
$R\mapsto \der_i(R)(x)$. As $\CC$-linear forms on $\CC [X]$,
the $(\der_i)_x$ are not linearly independent. There exists
a non-zero $p$-uple $(u_1,..,u_p)$ such that
$\sum_i u_i (\der_i)_x=0$. Since the $\partial_i$ are locally nilpotent
and commute pairwise, the derivation $\delta=u_1 (\der_1) +..+u_p (\der_p)$
is itself locally nilpotent. So $\delta$ generates the action $\varphi ^u$
defined by $\varphi ^u_t(y)=\varphi _{tu_1,..,tu_p}(y)$. Since
$\sum_i u_i (\der_i)_x=0$, $x$ is a fixed point of $\varphi ^u$
as can be seen via the exponential map. Therefore the stabilizer
of $x$ is not reduced to zero.
\qed
Let $\D=\{\der_1,..,\der_p\}$ be a commutative $p$-distribution on
$\CC[X]$. Since the exponential map defines a morphism of algebras,
the map:
$$
deg_D : \CC[X] \longrightarrow \mathbb{N} \cup \{-\infty \},\;
f \longmapsto deg_{t_1,..,t_p} \left \{ exp(t_1\der_1 +..+t_p\der_p)(f) \right
\}
$$
satisfies all the axioms of a degree function: This is the {\em degree
relative to $\D$}. By construction, the ring of invariants of $\D$ is
the set of regular functions of degree $\leq 0$. If $A$ is a domain,
we denote by $Fr(A)$ its fraction field. The following lemma is due to
Makar-Limanov ([M-L]).

\begin{lem} \label{Makar}
Let $A$ be a domain of characteristic zero. Let $\der$ be a non-zero
locally nilpotent derivation on $A$ and let $A^{\der}$ be its kernel.
Then $Fr(A)$ is isomorphic to $Fr(A^{\der})(t)$. In particular, for any
subfield $k$ of $Fr(A^{\der})$, the transcendence degrees satisfy the
relation:
$$
deg tr_k \{Fr(A^{\der})\}= deg tr_k \{Fr(A)\}-1
$$
\end{lem}
\dem Since $\der$ is non-zero locally nilpotent, there exists an element
$f$ of $A$ such that $\der(f)\not=0$ and $\der ^2 (f)=0$. So $g=\der(f)$
is invariant. It is then easy to check by induction on $p$ that every
element $P$ of $A$, of degree $p$ for $\der$, can be written
in a unique way as $g^p P = a_0 + ..+ a_p f^p$, where all the $a_i$
are invariant.
\qed
We end these recalls with the {\em factorial closedness} property, which
is essential for rings of invariants ([Da],[De]).

\begin{df}
Let $B$ a UFD and let $A$ be a subring of $B$. $A$ is factorially
closed in $B$ if every element $P$ of $B$ which divides a non-zero
element $Q$ of $A$ belongs to $A$.
\end{df}

\begin{lem}
Let $X$ be an affine variety such that $\CC[X]$ is a UFD. Let $\D$
be a commutative $p$-distribution on $X$. Then $\CC[X]^{\D}$
is factorially closed in $\CC[X]$.
\end{lem}
\dem Let $Q$ be a non-zero element of $\CC[X]^{\D}$, and let $P$ divide
$Q$ in $\CC[X]$. By considering the degree relative to $\D$, we get
$deg_{\D}(Q)=deg_{\D}(P) + deg_{\D}(Q/P)=0$. This implies
$deg_{\D}(P)=0$, and $P$ is invariant with respect to $\D$.
\qed

\section{Jacobian description of $p$-distributions}

Let $\varphi$ be an algebraic $\CP$-action on $\CN$, satisfying
the condition $(H)$. Let $\D$ be its commutative $p$-distribution,
and let $F$ be its quotient map. In this section we are going to
prove proposition \ref{Daigle}. The main idea is to construct a
system of rational coordinates for which calculations will be
simple. We obtain this system by adding some polynomials $s_i$ to
$f_1,..,f_{n-p}$. By analogy with $(\mathbb C,+)$-actions, we
denote them as "rational slices" (\cite{Da},\cite{D-F}). With
these coordinates, we show there exists an invariant fraction $E$
such that $[\D]=E\times J$, and there only remains to show that
$E$ is a polynomial.

\begin{df}
Let $\D$ be a commutative $p$-distribution on $\CX$. A diagonal
system of rational slices is a collection $\{s_1,..,s_p\}$ of
polynomials such that the matrix $({\partial}_i(s_j))$ is diagonal
and all its diagonal coefficients are non-zero invariant with
respect to $\D$.
\end{df}

\begin{lem} \label{Daigle5}
Every commutative $p$-distribution $\D$ satisfying the condition
$(H)$ admits a diagonal system of rational slices
$\{s_1,...,s_p\}$.
\end{lem}
\dem Let $\D_k$ be the commutative $(p-1)$-distribution
$\{\der_1,..,\der_{k-1},\der_{k+1},..,\der_p \}$, and let
$\CX^{\D_k}$ be its ring of invariants. By induction on lemma
\ref{Makar}, we get: $$ degtr_{\CC} Fr(\CX^{\D_k})\geq (n-p+1) $$
Since $\CX^{\D}$ is isomorphic to a polynomial ring in $(n-p)$
variables, $\der_k$ cannot be identically zero on $\CX^{\D_k}$.
For any $k$, there exists a polynomial $s_k$ such that
$\der_k(s_k)\not=0$, $\der_k ^2(s_k)=0$ and $\der_i(s_k)=0$ if
$i\not=k$. The collection $\{s_1,...,s_p\}$ is a diagonal system
of rational slices. \qed

\begin{lem}
Let $\D$ be a commutative $p$-distribution satisfying the
condition $(H)$. Let $\{s_1,..,s_p\}$ be a diagonal system of
rational slices. Then the map $G=(s_1,..,s_p,f_1,..,f_{n-p})$ is
dominating.
\end{lem}
\dem Let us show by absurd that $G$ is dominating. Assume that $G$ is not,
and let $Q$ be an element of $\CC[z_1,..,z_p,y_1,..,y_{n-p}]$ such that
$Q(G)=0$. We assume $Q$ to have minimal degree with respect to the variables
$z_1,..,z_p$. By derivation, we get for all $i$:
$$
\frac{\der Q}{\der  z_i}(G) \der_i(s_i)=\der_i(Q(G))=0
$$
Since $\der_i(s_i)\not=0$, this implies $\der Q /\der  z_i (G)=0$. By
minimality of the degree, we deduce that $\der Q /\der  z_i=0$ for all $i$.
So $Q$ belongs to $\CC[y_1,..,y_{n-p}]$. Therefore the $f_i$ are not
algebraically independent, and we obtain:
$$
degtr_{\CC} \CC(F) < n-p
$$
But $\CC[F]$ is the ring of invariants of $\D$. By induction with lemma
\ref{Makar}, we find that $degtr_{\CC} \CC(F)
\geq n-p$, hence a contradiction.
\qed

\begin{lem}
Let $\D$ be a commutative $p$-distribution satisfying $(H)$. Let
$\{s_1,..,s_p\}$ be a diagonal system of rational slices. Then
$\CX \subset \CC(f_1,..,f_{n-p})[s_1,..,s_p]$.
\end{lem}
\dem Let us show by induction on $r\geq 0$ that every polynomial
of degree $r$ with respect to $\D$ belongs to $\CC
(f_1,..,f_{n-p}) [s_1,..,s_p]$. For $r=0$, this is obvious because
every polynomial of degree zero is invariant, and belongs to $\CC
[f_1,..,f_{n-p}]$. Assume the property holds to the order $r$. Let
$R$ be a polynomial of degree $r+1$ with respect to $\D$. By
definition, the polynomials $\der_i(R)$ have all degree $\leq r$.
By induction, there exist some elements $P_i$ of
$\CC(y_1,..,y_{n-p})[z_1,..,z_p]$ such that $\der_i(R)=P_i(G)$ for
all $i$. Since $\D$ is commutative, we get for all $(i,j)$: $$
\frac{\der P_j}{\der z_i}(G)\der_i(s_i)=\der_i \circ \der_j (R)=
\der_j \circ \der_i (R)=\frac{\der P_i}{\der z_j}(G)\der_j(s_j) $$
By construction, there exists a non-zero polynomial $S_i$ in
$\CC[y_1,.., y_{n-p}]$ such that $\der_i(s_i)=S_i(F)$. Since $G$
is dominating, this yields for all $(i,j)$: $$ S_i\frac{\der
P_j}{\der z_i}=S_j\frac{\der P_i}{\der z_j} $$ The differential
1-form $\omega=\sum P_i/S_idz_i$ is polynomial in the variables
$z_i$. By the above equality, $\omega$ is closed with respect to
$z_i$. So $\omega$ is exact and there exists an element $P$ of
$\CC(y_1,..,y_{n-p})[z_1,...,z_p]$ such that $\omega=dP$.
Therefore $\der_i(R - P\circ G)=0$ for all $i$, and the function
$R - P\circ G$ is rational and invariant with respect to $\D$.
Since the ring of invariants of $\D$ is factorially closed, $R -
P\circ G$ belongs to $\CC(f_1,..,f_{n-p})$. So $R$ belongs to $\CC
(f_1,..,f_{n-p})[s_1,..,s_p]$, hence proving the induction. \qed
Following exactly the same argument, we can prove the equality: $$
\CX =\CC[f_1,..,f_{n-p}][s_1,..,s_p] $$ if the matrix
$(\der_i(s_j))$ is the identity. In this case $G$ is an algebraic
automorphism. In any case, the previous lemma asserts that $G$ is
always a birational automorphism of $\CN$.
\begin{lem} \label{Daigle3}
Let $\D$ be a commutative $p$-distribution satisfying $(H)$. Let
$\{s_1,...,s_p\}$ be a diagonal system of rational slices. Then
$\partial_1(s_1)..\partial_p(s_p)\times J= J(s_1,...,s_p)\times
[\D]$.
\end{lem}
\dem For any $p$-uple of polynomials $(R_1,...,R_p)$, there exist
some rational functions $P_i$ such that $R_i = P_i(G)$. On one hand,
we get by the chain rule:
$$
\begin{array}{ccl}
J(R_1,..,R_p)& = & \det(d(P_1,..,P_p,y_1,..,y_{n-p}))(G) \det(dG) \\ \\
 & = & \det((\der P_i /\der z_j))(G) J(s_1,..,s_p)
\end{array}
$$
On the other hand, we have the following relation:
$$
[\D](R_1,...,R_p)= \det((\der_i(R_j)))=\det (( \sum_k \der P_j / \der z_k (G)
\der_i(s_k)))
$$
Since the matrix $(\der_i (s_j))$ is diagonal, this yields:
$$
[\D](R_1,...,R_p)= \det ((\der P_i /\der z_j))(G)\der_1 (s_1)..\der_p (s_p)
$$
which implies the equality $\partial_1(s_1)..\partial_p(s_p) J(R_1,..,R_p)=
J(s_1,..,s_p)\times [\D](R_1,..,R_p)$.
\qed

\begin{lem} \label{Daigle4}
Let $\D$ be a commutative $p$-distribution satisfying the
condition $(H)$. Let $\{s_1,..,s_p\}$ be a diagonal system of
rational slices. Then $J(s_1,..,s_p)$ is invariant.
\end{lem}
\dem For simplicity, we denote by $J'$ the jacobian of every map
from $\CN$ to $\CN$. Since $\{s_1,..,s_p\}$ is a diagonal system
of rational slices, we get via the exponential map the relation
$s_i\circ \varphi= s_i + t_i\der_i(s_i)$, and this yields: $$
J'(s_1\circ\varphi,..,s_p\circ \varphi
,f_{1}\circ\varphi,..,f_{n-p}\circ\varphi)=
J'(s_1+t_1\der_1(s_1),..,s_p+t_p\der_p(s_p),f_{1},..,f_{n-p}) $$
Since every $\der_i(s_i)$ belongs to $\CC[F]$, we deduce: $$
J'(s_1\circ\varphi,..,s_p\circ\varphi,f_{1}\circ\varphi,..,f_{n-p}
\circ\varphi)= J'(s_1,..,s_p,f_{1},..,f_{n-p})=J(s_1,...,s_p) $$
Moreover we find by the chain rule: $$
J'(s_1\circ\varphi,..,s_p\circ\varphi,f_{1}\circ\varphi,..,f_{n-p}\circ\varphi)=
J'(s_1,..,s_p,f_{1},..,f_{n-p})(\varphi)\times J'(\varphi) $$
Since $\varphi$ is an automorphism of $\CN$ for any
$(t_1,...,t_p)$, the polynomial $J'(\varphi)$ never vanishes. So
it is non-zero constant. As $\varphi_{0,...,0}$ is the identity,
$J'(\varphi)\equiv 1$ and that implies: $$
J(s_1\circ\varphi,..,s_p\circ\varphi,f_{1}\circ\varphi,..,f_{n-p}\circ\varphi)=
J(s_1,..,s_p,f_{1},..,f_{n-p})(\varphi) $$ which leads to
$J(s_1,...,s_p)(\varphi)= J(s_1,...,s_p)$. Thus $J(s_1,...,s_p)$
is invariant. \qed {\bf{Proof of proposition \ref{Daigle}:}} Let
$\D$ be a commutative $p$-distribution satisfying the condition
$(H)$. By lemmas \ref{Daigle3} and \ref{Daigle4}, there exist two
non-zero invariant polynomials $E_1$ and $E_2$ such that: $$ E_1
\times [\D]=E_2 \times J $$ Since $\CC[F]$ is factorially closed
in $\CX$, we may assume that $E_1$ and $E_2$ have no common
factor. Let us show by absurd that $E_1$ is non-zero constant.
Assume that $E_1$ is not constant. By definition of $J$, $E_1$
divides all the coefficients of the $(n-p)$-form
$df_1\wedge..\wedge df_{n-p}$. So the hypersurface $V(E_1)$ is
contained in the singular set of $F$. But that contradicts a
result of Daigle (\cite{Da}), that asserts that $F$ is
non-singular in codimension 1. \qed

\section{Trivialisation of algebraic $\CP$-actions}

In this section, we are going to establish theorem \ref{Triv}. The
main idea is to refine a diagonal system of rational slices, in
order to get the coordinate functions of an algebraic automorphism
that conjugates $\varphi$ to the trivial action. \\ \ \\
{\bf{Proof of theorem \ref{Triv}:}} Let $\varphi$ be an algebraic
$\CP$-action on $\CN$ satisfying the condition $(H)$. Assume that
$E$ is constant and that the quotient map $F$ is quasi-fibered.
Let $\{s_1,..,s_p\}$ be a diagonal system of rational slices. Such
a system exists by lemma \ref{Daigle5}. By proposition
\ref{Daigle}, we have for any $(p-1)$-uple
$(R_1,..,R_{i-1},R_{i+1},..,R_p)$: $$
J(R_1,..,R_{i-1},s_i,R_{i+1},..,R_p)=[\D](R_1,..,R_{i-1},s_i,R_{i+1},..,R_p)/E
$$ Let $P_i$ be the polynomial of $\CC[t_1,..,t_{n-p}]$ such that
$\der_i(s_i)= P_i(F)$. Since $E$ is constant and $\der_k(s_i)=0$
if $k\not=i$, the previous equality yields: $$
J(R_1,..,R_{i-1},s_i,R_{i+1},..,R_p)\equiv 0 \;[P_i(F)] $$ If we
replace $R_k$ by all the polynomials $x_1,..,x_n$, we can see that
the coefficients of the differential form $ds_i \wedge df_1 \wedge
..\wedge df_{n-p}$ are all divisible by $P_i(F)$. By Daigle's
result (\cite{Da}), $F$ is non-singular in codimension 1. So the
coefficients of $df_1 \wedge ..\wedge df_{n-p}$ have no common
factor. Therefore $s_i$ satisfies the equation: $$ ds_i \wedge
\omega_F \equiv 0 \;[P_i(F)] $$ By the division lemma, there exist
some polynomials $A_i,S_i$ such that: $$ s_i = A_i(F) + P_i(F)S_i
$$ By an easy computation, we obtain that $(\der_i(S_j))$ is the
identity. By the remark following lemma \ref{Daigle3}, we have the
equality: $$ \CC[x_1,...,x_n]= \CC[f_1,..,f_{n-p}][S_1,..,S_p] $$
which implies that $G=(S_1,..,S_p,f_1,..,f_{n-p})$ is an algebraic
automorphism of $\CN$. Let $\varphi_0$ be the trivial action
generated by the commutative $p$-distribution $\{\der/\der
x_1,..,\der/\der x_p\}$. By using the exponential map, we find
that $G\circ \varphi = \varphi_0 \circ G$. So $\varphi$ is
trivial. \qed {\bf{Proof of corollary \ref{Triv2}:}} Let $\varphi$
be an algebraic $(\CC,+)$-action on $\CN$ satisfying $(H)$,
generated by the derivation $\der$. Assume that the quotient map
is quasi-fibered. Since $F$ is nonsingular in codimension 1, the
derivation $J$ is locally nilpotent and generates a
$(\CC,+)$-action $\varphi'$ such that $\cal{NL}(\varphi')$ has
codimension $\geq 2$. By theorem \ref{Triv}, $\varphi'$ is
trivial. Moreover via the automorphism of trivialisation, $\der$
is conjugate to $P(x_2,..,x_n)\der /\der x_1$, where $E=P(F)$ is
the factor of proposition \ref{Daigle}. \qed {\bf{Proof of
corollary \ref{Triv3}:}} Let $\varphi$ be an algebraic
$(\CC^{n-1},+)$-action on $\CN$, and assume that
${\cal{NL}}(\varphi)$ has codimension $\geq 2$. Then the factor
$E$ of proposition \ref{Daigle} is constant. Let us prove that
$\varphi$ is trivial. By theorem \ref{Triv}, we only have to show
that $\varphi$ satisfies the condition $(H)$ and that its quotient
map is quasi-fibered.

Let $f$ be a non-constant invariant polynomial of minimal
homogeneous degree on $\CX$. Then $f-\lambda$ is irreducible for
any $\lambda$. Indeed if $f-\lambda$ were reducible, all its
irreducible factors would be invariant by factorial closedness.
But that contradicts the minimality of the degree of $f$. Since
all the fibres of $f$ are irreducible, they are reduced and
connected. So $f$ is quasi-fibered, and there only remains to
prove that $f$ generates the ring of invariants of $\varphi$.

Let us show by induction on $r$ that any invariant polynomial $P$
of homogeneous degree $\leq r$ belongs to $\CC[F]$. This is
obvious for $r=0$. Assume this is true to the order $r$, and let
$P$ be an invariant polynomial of degree $\leq r+1$. Let $x$ be a
point in $\CN$ where $\varphi$ is free, and set $y=f(x)$. Since
$P$ is invariant, $P$ is constant on the orbit of $x$. Since this
orbit has dimension $(n-1)$ and that $f^{-1}(y)$ is irreducible,
this orbit is dense in $f^{-1}(y)$. So $P$ is constant on
$f^{-1}(y)$. By Hilbert's Nullstellensatz, there exists a
polynomial $Q$ such that $P=P(x)+ (f-y)Q$. The polynomial $Q$ is
invariant by factorial closedness and has degree $\leq r$. By
induction, $Q$ belongs to $\CC[F]$, and so does $P$, hence giving
the result. \qed

\section{A few examples}

We can show that the first assertion in theorem \ref{quasi3} is an
equivalence. More precisely, a primitive mapping $F$ is
quasi-fibered if and only if $\TF=0$. We will not prove it here,
but we would rather give two examples illustrating the necessity
of the conditions given in theorem \ref{Triv}. In both cases, the
module of relative exactness is not zero. Consider the locally
nilpotent derivation on $\CC[x,y,z]$: $$ \der _1 =x
\frac{\partial}{\partial y} - 2y\frac{\partial}{\partial z} $$ Its
ring of invariant is generated by $x$ and $xz+y^2$, and its
quotient map is defined by: $$ F_1: \CC ^3 \longrightarrow \CC ^2,
\quad (x,y,z)\longmapsto (x,xz+y^2) $$ It is easy to check that
$F_1$ is surjective and that $\overline{B(F_1)}= \{(u,v) \in \CC
^2, u=0\}$. So $F_1$ is not quasi-fibered because its fibres are
not 1-generically connected, and the action generated by $\der_1$
is not trivial. Second consider the locally nilpotent derivation
on $\CC [x,y,u,v]$: $$ \der _2 =u\frac{\partial}{\partial x} +
v\frac{\partial}{\partial y} $$ The polynomials $u,v,xv - yu$ are
invariant and generate the ring of invariants of $\der _2$. So the
corresponding action $\varphi_2$ satisfies the condition $(H)$,
and its quotient map is given by: $$ F_2 :\CC ^4 \longrightarrow
\CC ^3, \quad (x,y,u,v)\longmapsto (u,v,xv - yu) $$ By an easy
computation, we get that $B(F_2)$ is empty, $S(F_2)=V(x,y)$ and
$I(F_2)=\{(r,0,0), r \in \CC ^{*}\}$. So $F_2$ is not
quasi-fibered because its fibres are not 2-generically non-empty,
and $\varphi_2$ is not trivial.

\end{document}